\documentclass[leqno,12pt,a4paper]{amsart}

\usepackage{amsmath}
\usepackage{amsthm}
\usepackage{amssymb}
\usepackage{amscd}
\usepackage{bm}
\usepackage{graphicx,psfrag}

\theoremstyle{plain}
\newtheorem{Thm}{Theorem}[section]
\newtheorem{Lem}[Thm]{Lemma}

\newtheorem{Prop}[Thm]{Proposition}

\theoremstyle{definition}
\newtheorem{Def}[Thm]{Definition}
\newtheorem{Exa}[Thm]{Example}
\newtheorem{Rem}[Thm]{Remark}
\newtheorem{Hyp}[Thm]{Hypotheses}

\newcommand{\norm}[1]{\left|\!\left|#1\right|\!\right|}
\newcommand{\abs}[1]{\left|#1\right|}
\newcommand{\XX}{\sharp_{i\in\mathbb{Z}}X_i}
\newcommand{\deldel}{\frac{\partial}{\partial t}}

\title[Gluing an infinite number of instantons]{Gluing an infinite number of instantons}
\author[Masaki Tsukamoto]{Masaki Tsukamoto$^*$}

\date{\today}
\thanks{$^*$Supported by Grant-in-Aid for JSPS Fellows (19$\cdot$1530) from Japan Society for the
Promotion of Science}

\begin{document}

\subjclass[2000]{Primary 58D27; 
Secondary 53C07, 46T05}

\keywords{Yang-Mills gauge theory, gluing ASD connections, infinite energy, infinite dimensional moduli space, Donaldson's alternating method.} 

\maketitle
\begin{abstract}
This paper is one step toward infinite energy gauge theory and the geometry of infinite dimensional moduli spaces.
We generalize a gluing construction in the usual Yang-Mills gauge theory to an ``infinite energy'' situation.
We show that we can glue an infinite number of instantons, and that the resulting ASD connections have infinite energy in general. 
Moreover they have an infinite dimensional parameter space. 
Our construction is a generalization of Donaldson's ``alternating method''.
\end{abstract}

\section{Introduction}

Since the remarkable paper of S. K. Donaldson [D1] appeared in 1983, Yang-Mills gauge theory has been one of central themes in geometry and mathematical physics. There are plenty of papers on this subject, and many interesting theories have been developed.

But most of them study only \textit{finite energy} ASD connections and their \textit{finite dimensional} moduli spaces. 
 We know little about \textit{infinite energy} ASD connections and their moduli spaces. 
In fact I don't know even whether such connections exist in general situations. 
This paper is one step in this direction. 
We will show that we can construct a huge amount of infinite energy ASD connections by using a ``gluing construction''.

Gluing ASD connections is a well-known technique in Yang-Mills theory [D1], [D2], [DFK], [DK], [FU], [T1], [T2], but we usually study a gluing of two (or a finite number of) instantons.
In this paper we study a gluing of infinitely many instantons, and show that 
this new gluing technique produces the desired infinite energy ASD connections.
Moreover we show that the construction of ``gluing an infinite number of instantons" 
has the following interesting new feature;
To define the gluing of ASD connections, we need to introduce a ``gluing parameter space",
which is isomorphic to $SU(2)$ in the case of gluing two $SU(2)$ instantons.
In the case of gluing infinitely many $SU(2)$ instantons, the gluing parameter space becomes
the \textit{infinite product} of the copies of $SU(2)$.
This gigantic parameter space is the interesting new feature of our approach.

Although the above ``infinite dimensional gluing parameter space" is the source of the richness 
of the construction, it also causes several serious difficulties in analysis.
In particular, we cannot use the usual Fredholm deformation theory
because of infinite degrees of freedom of the parameter space.
We overcome this difficulty by adapting Donaldson's ``alternating method'' in [D2] to this situation. This decomposes the problem into an infinite number of small pieces where we can use the usual Fredholm deformation theory. 

I think that the technique of ``gluing an infinite number of instantons" 
probably can be applied to other equations
\footnote{After I posted a preprint version of this paper to the arXiv in 2005, 
I was informed that there is a paper of M. Macr\`{i}, M. Nolasco and T. Ricciardi [MNR] which 
studies an ``infinite energy" situation for selfdual vortices in $\mathbb{R}^2$.
From [MNR] I learned that there is a paper of S. Angenent [A] 
which develops a gluing technique for infinitely many solutions of an elliptic PDE 
of the form $\Delta u + f(x, u) = 0$ in $\mathbb{R}^n$.
[MNR] adapts this gluing technique to an elliptic equation with singular sources.
I think that our alternating method is different from the analytic technique in [A].} , e.g., the Seiberg-Witten equation, pseudo-holomorphic curves, etc.
In fact the classical Mittag-Leffler theorem in complex analysis can be considered as a variant of 
this gluing technique. 
Singularities of a meromorphic function are instanton-like objects, and
the Mittag-Leffler theorem constructs a meromorphic function which has 
an infinite number of singularities given in the complex plane.
In other words the result of this paper is a Yang-Mills analogue of the Mittag-Leffler theorem.

The organization of the paper is as follows.
In Section 2 we present notation and state  main theorems. We establish basic estimates in Section 3.
The details of the gluing construction are given in Section 4. We discuss a moduli problem in Section 5.

\vspace{2mm}

\textbf{Acknowledgement.}\hspace{0.7mm}
 I wish to thank Professor Kenji Fukaya for his suggestion and help during the course of this work.


\section{Main Results}
\subsection{Brief review of Yang-Mills theory}\label{subsection: review of Yang-Mills}
To begin with, we review some basic notions of Yang-Mills theory. For the details, see [DK] or [FU].
Let $G$ be a Lie group $SU(2)$ or $SO(3)$.  Let $C(G)$ be the center of $G$,
 i.e., $C(SU(2)) = \{\pm 1\}$ and $C(SO(3)) = 1$.
Suppose that $X$ is a closed oriented Riemannian 4-manifold and that 
$E \rightarrow X$ is a principal $G$ bundle over $X$. 
Let $\mathrm{ad}E$ be the Lie algebra bundle associated with $E$.
A connection $A$ on $E$ is said to be irreducible if the isotropy group of $A$ in the gauge group 
is equal to $C(G)$. ($C(G)$ acts trivially on all connections.)

A connection $A$ is said to be anti-self-dual, or ASD, if the  self-dual part of its curvature is zero:
\[ F^+_A = 0.\]
If $A$ is an ASD connection, then we have the Atiyah-Hitchin-Singer complex for $A$:
\[ 0 \rightarrow \Omega^0(\mathrm{ad}E) \xrightarrow{d_A} \Omega^1(\mathrm{ad}E) \xrightarrow{d^+_A} \Omega^+(\mathrm{ad}E) \rightarrow 0,
\]
where $\Omega^+$ is the space of self-dual 2-forms and $d^+_A$ denotes the self-dual part of $d_A$.
An ASD connection $A$ is said to be acyclic if the above complex is exact (acyclic).
The Atiyah-Hitchin-Singer complex is exact if and only if 
\[d^*_A + d^+_A : \Omega^1(\mathrm{ad}E) \rightarrow 
\Omega^0(\mathrm{ad}E) \oplus \Omega^+(\mathrm{ad}E) \]
is isomorphism. Here $d^*_A$ is the formal adjoint of $d_A : \Omega^0(\mathrm{ad}E) \rightarrow \Omega^1(\mathrm{ad}E)$. 
$A$ is said to be regular if $A$ is irreducible and $d^*_A + d^+_A$ is surjective.
\begin{Exa}
Let $X$ be the 4-sphere $S^4$ with the usual round metric.
All non-flat $G$ ASD connections on $S^4$ are regular.
This follows from the Weitzenb\"{o}ck formula. For example, see [FU, Chapter 6].
\end{Exa}
\begin{Exa}
$K3$ surface has a non-flat, irreducible, acyclic $SO(3)$ ASD connection. See [DK, Section 9.1.3].
\end{Exa}

\subsection{Gluing construction}\label{subsection: gluing construction}
Next we introduce the following convenient definition.
 \begin{Def}\label{def: gluing data}
A 5-tuple $(X,x^L,x^R,E,A)$ is called a gluing data if it satisfies the following (i)$\sim$(v):
     \begin{enumerate}
       \item $X$ is a closed oriented Riemannian 4-manifold.
       \item $x^L$ and $x^R$ are two distinct points in $X$.
       \item $E$ is a principal $G$ bundle over $X$.
       \item $A$ is an irreducible acyclic ASD connection on $E$. 
(We can replace this condition with the weaker condition that $A$ is a regular ASD connection. 
See the remark at the end of Section 5.) 
       \item The metric on $X$ is flat in some neighborhoods of $x^L$ and $x^R$. (See the remark below.)
    \end{enumerate}

\end{Def}
\begin{Rem}
The condition (v) in the Definition \ref{def: gluing data} is artificial and just for simplicity.
We can construct our gluing argument without this condition; see [D2, (4.48) Proposition].
Actually, this condition is allowable in the following sense.
Let $g$ be a Riemannian metric on $X$ which does \textit{not} satisfy the condition (v).
Let $A$ be an irreducible acyclic (resp. regular) $g$-ASD connection on $E$.
Then small perturbation gives a new metric $g'$ satisfying the condition (v) and 
an irreducible acyclic (resp. regular) $g'$-ASD connection $A'$.
This fact follows from the usual implicit function theorem.
\end{Rem}
\begin{Def}\label{def: finite type}
   Let $(X_i,x^L_i,x^R_i,E_i,A_i)$, $(i \in \mathbb{Z})$, be a sequence of gluing data indexed by integers.
 This sequence is said to be \textit{of finite type} if it satisfies the following: 
\[ \{(X_i,x^L_i,x^R_i,E_i,A_i) |i \in \mathbb{Z}\} \text{ is a  finite set}. \]
(Of course, the sequence itself is an \textit{infinite} sequence.) 
This condition assures the uniformity of analysis and geometry on each data.
\end{Def}
We will define a connected sum of gluing data. This construction is a generalization of the ``conformal connected sum'' construction in [D2, pp. 306-307] and [DK, Section 7.2.1].

\begin{figure}[h]
      \psfrag{1}{$X$}
      \psfrag{2}{$x^L$}
      \psfrag{3}{$x^R$}
      \psfrag{4}{$U$} 
      \psfrag{5}{$\Omega^L$}
      \psfrag{6}{$\Omega^R$}
      \psfrag{7}{$L^-$}
      \psfrag{8}{$L^+$}
      \psfrag{9}{$R^+$}
      \psfrag{10}{$R^-$}
      \psfrag{11}{$U_i$}
      \psfrag{12}{$U_{i+1}$}
      \psfrag{13}{$R^+_i = L^-_{i+1}$}
      \psfrag{14}{$R^-_i = L^+_{i+1}$}
      \psfrag{15}{$\Omega^R_i = \Omega^L_{i+1}$}
      \psfrag{16}{$x^R$}
      \psfrag{17}{$R^-$}
      \psfrag{18}{$R^+$}
      \psfrag{19}{$\Omega^R$}
      \psfrag{20}{$ k N^{-1} \sqrt{\lambda}$}
      \psfrag{21}{$ N^{-1} \sqrt{\lambda}$}
      \psfrag{22}{$N \sqrt{\lambda}$}
      \psfrag{23}{$ k^{-1} N \sqrt{\lambda}$}
      \psfrag{24}{$\XX$}     

     \begin{flushleft}
       \includegraphics[scale = 1]{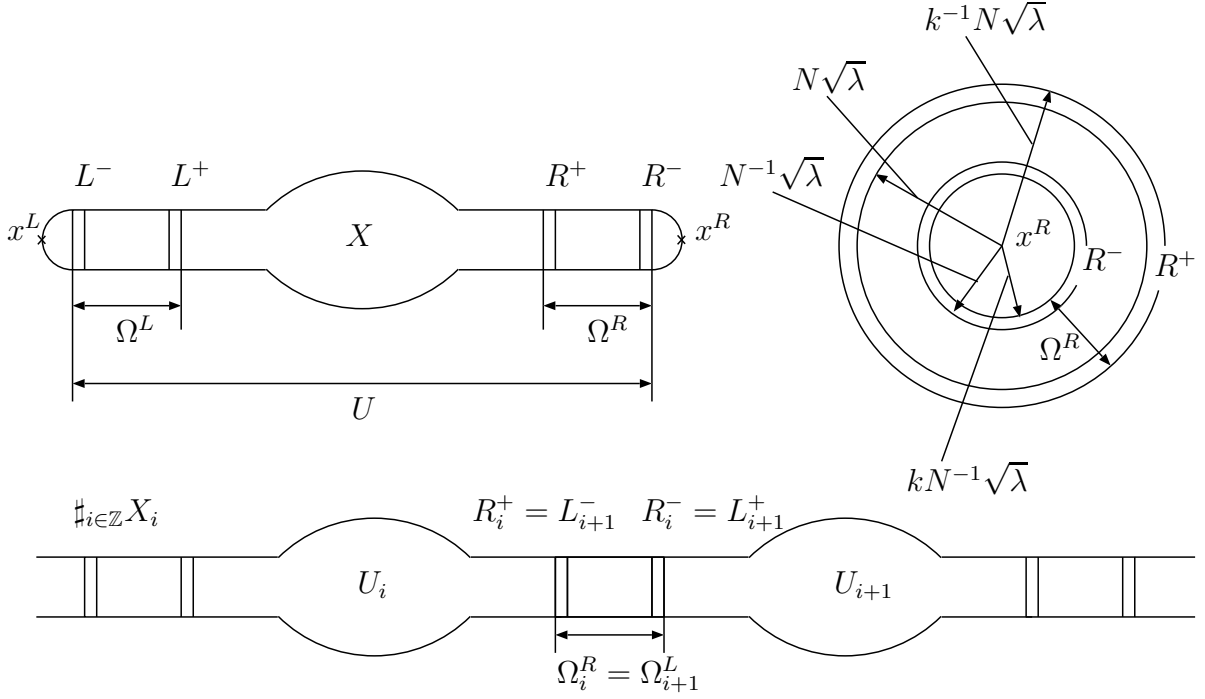}
     \end{flushleft}  
     \caption{A gluing data and a connected sum.}  \label{fig: gluing data}
     \end{figure}

 Let $(X,x^L,x^R,E,A)$ be a gluing data. By the condition (v) in Definition \ref{def: gluing data}, there are an oriented Euclidean coordinate $\xi$ centered on $x^R$ and 
an oriented Euclidean coordinate $\eta$ centered on $x^L$.
We set $k := 0.9$, and 
we define the annular region $\Omega^R$ and the ``shells'' $R^-$ and $R^+$ in the neighborhood of $x^R$  by
\begin{equation*}
 \begin{split} 
  \Omega^R &:= \{ \xi \,|\, kN^{-1}\sqrt{\lambda}< |\xi| < k^{-1}N\sqrt{\lambda}\} ,\\
  R^- &:= \{\xi \,|\, kN^{-1}\sqrt{\lambda}< |\xi| < N^{-1}\sqrt{\lambda}\} ,\\
  R^+ &:= \{\xi \,|\, N\sqrt{\lambda}< |\xi| < k^{-1}N\sqrt{\lambda}\} .
 \end{split}
\end{equation*}
Here $N$ and $\lambda$ are positive parameters.
We define $\Omega^L$, $L^-$ and $L^+$ in the same way by using the Euclidean coordinate $\eta$ centered on $x^L$. 
So $L^-$ is the inner ``shell" of $\Omega^L$ and $L^+$ is the outer ``shell";
 see Figure \ref{fig: gluing data}. 
(Figure \ref{fig: gluing data} is a variant of [D2, Diagram (4.13)].)
We define the open set $U$ by 
\[ U := X \setminus (\Bar{B}(x^R, kN^{-1}\sqrt{\lambda}\,)\cup \Bar{B}(x^L,kN^{-1}\sqrt{\lambda}\,)). \]  
Here $\Bar{B}(x^R, kN^{-1}\sqrt{\lambda}\,)$ and $\Bar{B}(x^L,kN^{-1}\sqrt{\lambda}\,)$ 
are the balls $|\xi|\leq kN^{-1}\sqrt{\lambda}$ and $|\eta|\leq kN^{-1}\sqrt{\lambda}$.  

In the later gluing argument, ``error terms" will be supported in the shells $R^{\pm}$ and $L^{\pm}$.
The volumes of these shells are determined by the parameters $\lambda$ and $N$.
We will choose $\lambda$ small and $N$ large. 
The value of $N$ will be fixed in the end of Section 3. 
We choose $\lambda$ so small that 
\begin{equation}
\lambda N^{100} \ll 1. \label{lambda is small}
\end{equation}
Here $``100"$ does not have a particular meaning. The point is that $\lambda$ is very small.
(In the above, we set $k = 0.9$. This value $``0.9"$ does not have a particular meaning, either.
The point is that $k$ is smaller than 1.)

 Let $(X_i,x^L_i,x^R_i,E_i,A_i)_{i\in \mathbb{Z}}$ be a sequence of gluing data of finite type. 
We have $\Omega^R_i$, $R^-_i$, $R^+_i$, and $\Omega^L_i$, $L^-_i$, $L^+_i$, and $U_i$ for each $i\in\mathbb{Z}$. 
 We define the connected sum
\[\sharp_{i\in\mathbb{Z}}X_i :=  \cdots \sharp X_{-1} \sharp X_0 \sharp X_1 \sharp
\cdots.
\]
by identifying $\Omega^R_i$ and $\Omega^L_{i+1}$. 
The identification map between $\Omega^R_i$ and $\Omega^L_{i+1}$ is defined by
\begin{equation}
 \eta = \lambda \bar{\xi}/|\xi|^2.  \label{conformal map}
\end{equation}
Here $\xi$ is an oriented Euclidean coordinate centered on $x^R_i$, 
and $\eta$ is an oriented Euclidean coordinate centered on $x^L_{i+1}$.
$\xi$ $\mapsto$ $\bar{\xi}$ is any reflection.
$R^-_i$ and $R^+_i$ are identified with $L^+_{i+1}$ and $L^-_{i+1}$, respectively. 
The connected sum space $\XX$ is covered by the open sets $U_i$:
\[ \sharp_{i\in\mathbb{Z}}X_i = \bigcup_{i\in \mathbb{Z}}U_i. \]
Since the identification map (\ref{conformal map}) is conformal,
$\sharp_{i\in\mathbb{Z}}X_i$ has a natural conformal structure which is compatible with the metrics of the given manifolds $X_i$. 
(Note that ASD condition is conformally invariant.) 

Next we will define a gluing of the principal $G$ bundles $E_i$.
\begin{Def}
We define the \textit{gluing parameter space} GlP 
and the \textit{effective gluing parameter space} EGlP by
\begin{equation*}
 \begin{split}
  \mathrm{GlP} &:= \prod_{i\in\mathbb{Z}}\mathrm{Hom}_{G}((E_i)_{x^R_i},(E_{i+1})_{x^L_{i+1}})
  \; \cong G^{\mathbb{Z}},\\
  \mathrm{EGlP} &:= \prod_{i\in\mathbb{Z}}(\mathrm{Hom}_{G}((E_i)_{x^R_i},(E_{i+1})_{x^L_{i+1}})/C(G))
  \; \cong (G/C(G))^{\mathbb{Z}}.
 \end{split}
\end{equation*}
Here
$\mathrm{Hom}_{G}((E_i)_{x^R_i},(E_{i+1})_{x^L_{i+1}})$ is the set of $G$-isomorphisms between the fibres $(E_i)_{x^R_i}$ and $(E_{i+1})_{x^L_{i+1}}$, and it is isomorphic to $G$.
\end{Def}
If $x\in X_i$ is a point in a small neighborhood of $x^R_i$, 
then the fibre $(E_i)_{x}$ can be identified with $(E_i)_{x^R_i}$ 
by using the parallel transport along the radial line from $x^R_i$ to $x$.
This ``exponential gauge" centered on $x_i^R$ defines a local trivialization of $E_i$ over $\Omega^R_i$. 
In the same way 
we define a trivialization of $E_i$ over $\Omega^L_i$ by the exponential gauge centered on $x_i^L$. 
(For the detail of exponential gauge, see [FU, Chapter 9] or [DK, Section 2.3.1].)
Let $\rho = (\rho_i)_{i\in \mathbb{Z}}\in \mathrm{GlP}$ be a gluing parameter.
$\rho_i$ is an isomorphism between the fibres $(E_i)_{x^R_i}$ and $(E_{i+1})_{x^L_{i+1}}$.
Since we have the identification $\Omega^R_i \cong \Omega^L_{i+1}$ in (\ref{conformal map}) and 
the bundle trivializations $E_i|_{\Omega^R_i} \cong \Omega^R_i \times (E_i)_{x^R_i}$ and
$E_{i+1}|_{\Omega^L_{i+1}} \cong \Omega^L_{i+1} \times (E_{i+1})_{x^L_{i+1}}$,
$\rho_i$ defines an identification map between $E_i|_{\Omega^R_i}$ and $E_{i+1}|_{\Omega^L_{i+1}}$
covering the base space identification (\ref{conformal map}).
(See the diagram below.)
We define a principal $G$ bundle $\bm{E}(\rho)$ over $\sharp_{i\in\mathbb{Z}}X_i$ by using this identification map. 
\begin{equation}\label{bundle gluing}
\begin{CD}
E_i|_{\Omega^R_i} @>{\text{exp. gauge}}>>  \Omega^R_i \times (E_i)_{x^R_i}\\
@VVV @VV{\rho_i}V \\
E_{i+1}|_{\Omega^L_{i+1}} @>{\text{exp. gauge}}>>  \Omega^L_{i+1} \times (E_{i+1})_{x^L_{i+1}}
\end{CD}
\end{equation}


The gluing construction in Section 4 will give an ASD connection $\bm{A}(\rho)$ on $\bm{E}(\rho)$ which is ``close'' to $A_i$ over each $U_i$:
\begin{Thm}\label{thm: construction}
For any positive number $\varepsilon >0$ and $p \ge 1$, 
we can choose the parameters $N$ and $\lambda$ so that 
there exists an ASD connection $\bm{A}(\rho)$ on $\bm{E}(\rho)$ for all $\rho \in \mathrm{GlP}$  
such that 
\[ |\!|\bm{A}(\rho) - A_i|\!|_{L^p(U_i)} < \varepsilon \quad \text{for all}\ i \in \mathbb{Z}.  \]  
(Since there is a natural identification between $\bm{E}(\rho)$ and $E_i$ over $U_i$,
$\bm{A}(\rho) - A_i$ is well defined over $U_i$.)
Moreover, if there are infinitely many non-flat ASD connections $A_i$ in the given sequence of gluing data, 
then the ASD connection $\bm{A}(\rho)$ has infinite $L^2$-energy:
\[\norm{F(\bm{A}(\rho))}_{L^2} = + \infty. \]
\end{Thm}
In the above, the norm $\norm{\cdot}_{L^p(U_i)}$ is defined by using the Riemannian metric on $X_i$.
Since the $L^2$-norm of a 2-form is conformally invariant, the $L^2$-norm of the curvature $\norm{F(\bm{A}(\rho))}_{L^2}$
is well-defined. ($\sharp_{i\in\mathbb{Z}}X_i$ has the natural conformal structure.)
We will prove a more precise result in Section 4.

In Section 5 we will discuss a moduli problem. Our main result is the following. 
(Here the parameters $N$ and $\lambda$ are chosen appropriately.)
\begin{Thm}\label{thm: injectivity}
Let $\rho$ and $\rho'$ be two gluing parameters. 
Then $\bm{A}(\rho)$ is gauge equivalent to $\bm{A}(\rho')$ if and only if $[\rho] = [\rho']$ in $\mathrm{EGlP}$.
Here the gauge equivalence is defined by usual smooth bundle maps between the principal $G$ bundles $\bm{E}(\rho)$ and $\bm{E}(\rho')$.
\end{Thm} 
The effective gluing parameter space $\mathrm{EGlP} = (G/C(G))^{\mathbb{Z}}$ has infinite degrees of freedom. Hence we can consider that our ASD connections $\bm{A}(\rho)$ have infinite dimensional parameter space.

\subsection{Outline of the alternating method}\label{subsection: outline}
We will construct the ASD connection $\bm{A}(\rho)$ by applying
``Donaldson's alternating method" in [D2]. Here we explain the basic idea of it.
The alternating method is an iterative process for solving the ASD equation.
First we construct the initial approximate solution $\bm{A}^{(0)}(\rho)$ on $\bm{E}(\rho)$ 
by using cut-off functions.
Of course, this $\bm{A}^{(0)}(\rho)$ is not ASD in general, i.e., $F^{+}(\bm{A}^{(0)}(\rho)) \neq 0$.
The alternating method iteratively improves this $\bm{A}^{(0)}(\rho)$ as follows;
\begin{description}
\item[Step 0]
To cancel the ``error term" $F^{+}(\bm{A}^{(0)}(\rho))$,
we slightly perturb $\bm{A}^{(0)}(\rho)$ over each $X_{2i}$
and construct $\bm{A}^{(1)}(\rho)$ on $\bm{E}(\rho)$ 
by using these perturbations and cut-off functions.
The new error term $F^{+}(\bm{A}^{(1)}(\rho))$ becomes smaller than 
the old error term $F^{+}(\bm{A}^{(0)}(\rho))$.
\item[Step 1]
To cancel the error term $F^{+}(\bm{A}^{(1)}(\rho))$,
we slightly perturb $\bm{A}^{(1)}(\rho)$ over each $X_{2i+1}$
and construct $\bm{A}^{(2)}(\rho)$ on $\bm{E}(\rho)$ 
by using these perturbations and cut-off functions.
The new error term $F^{+}(\bm{A}^{(2)}(\rho))$ becomes smaller than 
the old error term $F^{+}(\bm{A}^{(1)}(\rho))$.
\end{description}
Step 2 is the same as Step 0, and Step 3 is the same as Step 1. In general, 
\begin{description}
\item[Step $(2n)$]
We slightly perturb $\bm{A}^{(2n)}(\rho)$ over each ``even" manifold $X_{2i}$
and construct $\bm{A}^{(2n +1)}(\rho)$.
\item[Step $(2n+1)$]
We slightly perturb $\bm{A}^{(2n+1)}(\rho)$ over each ``odd" manifold $X_{2i+1}$
and construct $\bm{A}^{(2n+2)}(\rho)$.
\end{description}
The error terms $F^{+}(\bm{A}^{(n)}(\rho))$ go to $0$ as $n$ goes to infinity,
and the sequence of the approximate solutions $\bm{A}^{(n)}(\rho)$ converges to 
the desired ASD connection $\bm{A}(\rho)$.
This is the outline of the alternating method in Section 4.

\section{Main Estimates}
In this section we establish basic estimates for the proofs of Theorem \ref{thm: construction} and 
\ref{thm: injectivity}. 
These are essentially the reproduction of the estimates in [D2, (4.22) Lemma and (4.24) Remarks]. 
We follow the arguments of [D2, IV (iii)].

Fix $p >6$ and let $(X, x^L, x^R, E, A)$ be a gluing data. 
We have open sets $L^{\pm}$, $R^{\pm}$ and $U$ on $X$. 
Let us consider the following map over $X$:
\begin{equation}\label{map: perturb}
 \begin{split}
 L^p_1(\Omega^1 (\mathrm{ad}E)) &\times L^{2p}(\Omega^1 (\mathrm{ad}E))
 \rightarrow L^p(\Omega^0 (\mathrm{ad}E) \oplus  \Omega^+ (\mathrm{ad}E)),  \\
 &( b , a ) \mapsto D_A b + [a \wedge b]^+ + (b \wedge b)^+. 
 \end{split}
\end{equation}
Here $D_A$ is the differential operator $d^*_A + d^+_A : \Omega^1(\mathrm{ad}E)
 \rightarrow \Omega^0 (\mathrm{ad}E) \oplus \Omega^+(\mathrm{ad}E)$.
Since we suppose that $A$ is acyclic in Definition \ref{def: gluing data}, 
$D_A$ becomes an invertible map from $L^p_1$ to $L^p$.
Hence if $a \in \Omega^1(\mathrm{ad}E)$ is small in $L^{2p}$ 
and $\sigma \in \Omega^+(\mathrm{ad}E)$ is small in $L^{p}$,
then we can solve the following equation for $b\in \Omega^1(\mathrm{ad}E)$:
\begin{equation} \label{eq: perturb}
D_A b + [a \wedge b]^+ +(b \wedge b)^+ = - \sigma, 
\end{equation}\label{apriori estimate for b}
and we get a solution $b$ which satisfies
\begin{equation}
|\!|b|\!|_{L^p_1(X)} \leq  \mathrm{const} \cdot  |\!|\sigma|\!|_{L^p(X)}. \label{implicit}
\end{equation}
(This is because $b = 0$ is a solution to (\ref{eq: perturb}) if $\sigma = 0$.)
 
Suppose that $a$ is a smooth form defined on $U$ and that the self-dual part of the curvature
$\sigma := F^+(A + a)$ is supported in $L^+ \cup R^+$. ($\sigma$ is also defined on $U$.) 
We extend $a$ and $\sigma$ to $X$ by zero, and set
\[ \delta := |\!|\sigma|\!|_{L^\infty (X)} .\] 
The volume of $(L^+\cup R^+)$ which contains the support of $\sigma$ is bounded by 
$\mathrm{const}\cdot  N^4 \lambda^2$.
Hence if $a$ is small in $L^{2p}$ and $\delta$ is bounded by some constant, 
then we have a small solution $b$ to (\ref{eq: perturb}) for sufficiently small $\lambda$.
In the later gluing argument $\sigma$ will be the ``error term" and 
$b$ will be the ``perturbation" described in Section \ref{subsection: outline}.

Let us introduce a cut-off function $\psi$ which satisfies $\mathrm{supp}(d\psi) \subset L^-\cup R^-$ and
 $\psi =0$ on $X\setminus U$ and $\psi = 1$ on $U\setminus(L^-\cup R^-)$. 
We can choose $\psi$ such that
\begin{equation}
|d\psi| \leq  \frac{20 N}{\sqrt{\lambda}}. \label{cutoff}
\end{equation}
(Note that $20N/\sqrt{\lambda}$ is a large number because we choose $\lambda$ very small.)
We define the ``new error term" $\tau$ by $\tau := F^+(A + a + \psi b)$ on $U$. 
\begin{eqnarray*}
\tau &=& F^+(A + a) + d^+_A(\psi b) + \psi [a\wedge b]^+ + \psi^2(b\wedge b)^+\\
     &=& \sigma + (d\psi \wedge b)^+ + \psi (d^+_A b + [a\wedge b]^+ + (b\wedge b)^+) + \psi(\psi - 1)(b\wedge b)^+\\
     &=& (1-\psi)\sigma + (d\psi \wedge b)^+ + \psi(\psi - 1)(b\wedge b)^+\\
     &=& (d\psi \wedge b)^+ + \psi(\psi - 1)(b\wedge b)^+, \quad 
\text(\mathrm{supp}(\sigma)\cap \mathrm{supp}(1-\psi) = \emptyset).
\end{eqnarray*}
Hence $\tau$ is supported in $L^- \cup R^-$, and we have 
\begin{equation}\label{estimate of tau by b}
|\!|\tau|\!|_{L^{\infty}} \leq \frac{20N}{\sqrt{\lambda}}\ |\!| b|_{L^-\cup R^-}|\!|_{L^\infty} + |\!| b|_{L^-\cup R^-}|\!|_{L^\infty}^2.
\end{equation}
Note that the ``old" error term $\sigma$ is supported in $L^+ \cup R^+$ 
and that the ``new" error term $\tau$ is supported in $L^- \cup R^-$.
\begin{Prop}\label{prop: estimate for construction}
For any $\kappa > 0$, there is $\eta > 0$ such that if $\delta \leq \kappa$, $|\!|a|\!|_{L^{2p}} \leq \eta$ and $\lambda \ll 1$, 
then we have the following:
\begin{eqnarray*}
&\mathrm{(i)}&
 |\!| b |\!|_{L^{2p}}\leq C_1(N)\cdot \lambda^{(p+2)/{2p}} \cdot \delta.\\
&\mathrm{(ii)}&
 |\!| \tau |\!|_{L^\infty} \leq C_2 \cdot K_N \cdot N^4 \cdot \delta.
\hspace{6mm} \text{Here we put }\hspace{1mm} K_N := N/(N-N^{-1})^3.\\
&\mathrm{(iii)}&
 |\!| \psi \cdot  b |\!|_{L^p_1} \leq C_3(N) \cdot \lambda^{2/p} \cdot \delta.\\
&\mathrm{(iv)}&
 |\!|F(A + a + \psi \cdot b) - F(A + a)|\!|_{L^2(U)} \leq C_4(N) \cdot \lambda \cdot \delta.
\end{eqnarray*}
Here $C_1(N)$, $C_3(N)$, $C_4(N)$ are constants depending on N,
and $C_2$ is a constant independent of $N$.
(The fact that $C_2$ is independent of $N$ is essential to the later argument.)
\end{Prop}
\begin{proof}   
The proof follows [D2, (4.22) Lemma]. In our case, some care must be paid to interactions between $L^{\pm}$ and $R^{\mp}$.  

\noindent            
(i) This follows from (\ref{implicit}) and the Sobolev embedding $L^{4p/(p+2)}_1 \hookrightarrow L^{2p}.$ (Note that this Sobolev embedding \textit{does not} hold at p = $\infty.$)
\begin{eqnarray*}
\norm{b}_{L^{2p}} &\leq& \mathrm{const} \cdot |\!| b |\!|_{L^{4p/(p+2)}_1}\\
                   &\leq& \mathrm{const} \cdot |\!| \sigma |\!|_{L^{4p/(p+2)}}\\
                   &\leq& C_1(N) \cdot \lambda^{(p+2)/2p}\cdot \delta.
\end{eqnarray*}
\noindent
(ii) The estimate of the $L^{\infty}$ norm of $\tau$ can be derived from the $L^{\infty}$ norm of $b$ 
over $L^- \cup R^-$. (See (\ref{estimate of tau by b}).)
First we will estimate a ``leading term" of $b$.
Let $\tilde{b}$ be the solution to 
\[ D_A \tilde{b} = - \sigma. \]
Since $D_A = d^*_A \oplus d^+_A$ is invertible, 
we have the Green kernel $G(x,y)$ for $D_A$:
\[\tilde{b}(x) = -\int_{X} G(x,y)\cdot \sigma (y)\, dy.\] 
$G(x,y)$ has a singularity along the diagonal such that  
\[ |G(x,y)| \leq \mathrm{const} \cdot \frac{1}{d(x,y)^3}. \]
Here $d(x,y)$ is the distance between $x$ and $y$. 
This is because $D_A$ is a \textit{first order} elliptic differential operator.
(For the details, see [D2, (4.22) Lemma].) 
In the complement of a neighborhood of the diagonal, $G(x,y)$ is bounded.

 Let us consider the case $x \in R^-$. Since $\sigma$ is supported in $L^+ \cup R^+$,
\[\tilde{b}(x) = - \int_{L^+} G(x,y)\cdot \sigma(y) \, dy - \int_{R^+} G(x,y)\cdot \sigma(y)\, dy.\]
 Because the distance between $R^-$ and $L^+$ is  sufficiently large, the first term can be easily estimated:
\[
\abs{\int_{L^+} G(x,y)\cdot \sigma(y) \, dy}\leq \mathrm{const}\cdot \mathrm{vol}(L^+)\cdot \delta 
                           \leq \mathrm{const}\cdot N^4\cdot \lambda^2\cdot \delta.
\]
This implies that the interaction between $R^-$ and $L^+$ is comparatively small.
The second term is more delicate. 
Since $d(R^+, R^-) = N\sqrt{\lambda}-N^{-1} \sqrt{\lambda}$, we have
\begin{eqnarray*}
\abs{\int_{R^+} G(x,y)\cdot \sigma(y) \, dy} &\leq& \mathrm{const}\cdot \mathrm{vol}(R^+)\cdot \frac{\delta}{(N\sqrt{\lambda} - N^{-1} \sqrt{\lambda})^3}\\
&\leq& \mathrm{const}\cdot K_N \cdot N^3 \cdot \sqrt{\lambda} \cdot \delta ,
\quad (K_N = \frac{N}{(N-N^{-1})^3}).
\end{eqnarray*}
Using (\ref{lambda is small}) in Section 2.2, we get
\[ |\tilde{b}(x)|\leq \mathrm{const}\cdot K_N \cdot N^3 \cdot \sqrt{\lambda}\cdot \delta.\]
We have the same estimate in the case $x \in L^-$.

Next we will estimate the ``lower term". 
Set $b=\tilde{b}+\beta$. Then $\beta$ satisfies
\[D_A\beta + [a\wedge b]^+ +(b\wedge b)^+ = 0.\]
From the Sobolev embeddings,
\begin{eqnarray*}
|\!| \beta |\!|_{L^{\infty}} &\leq& \mathrm{const}\cdot \norm{D_A\beta}_{L^6},
 \quad (L^6_1 \hookrightarrow L^{\infty}),\\
&\leq& \mathrm{const} \cdot(|\!|a|\!|_{L^{12}}\cdot|\!|b|\!|_{L^{12}} + |\!|b|\!|^2_{L^{12}})\\
&\leq& \mathrm{const} \cdot(|\!|a|\!|_{L^{2p}}\cdot |\!|b|\!|_{L^3_1} + |\!|b|\!|^2_{L^3_1}),
 \quad (L^3_1\hookrightarrow L^{12}),\\
&\leq& \mathrm{const}\cdot(\eta \cdot |\!|b|\!|_{L^3_1} + |\!|b|\!|^2_{L^3_1}).
\end{eqnarray*}
From (\ref{implicit}) we have
\[|\!|b|\!|_{L^3_1}\leq \mathrm{const}\cdot |\!|\sigma|\!|_{L^3} \leq \mathrm{const}\cdot N^{4/3}\cdot \lambda^{2/3}\cdot \delta.\]
Hence
\[|\!| \beta |\!|_{L^{\infty}} \leq C(N)\cdot \lambda^{2/3}\cdot \delta ,
 \quad(C(N) \text{ is a constant depending on N}).
\]
Since $C(N)\cdot \lambda^{2/3} \ll \sqrt{\lambda}$, $\tilde{b}$ 
is the leading term of $b$ over $R^- \cup L^-$:
\begin{equation}\label{L^infty of b}
\norm{ b|_{R^-\cup L^-}}_{L^{\infty}} \leq \mathrm{const}\cdot K_N\cdot N^3\cdot \sqrt{\lambda}\cdot \delta.
\end{equation}
Substituting this into (\ref{estimate of tau by b}), we conclude that
\[|\!|\tau|\!|_{L^{\infty}}\leq \frac{20N}{\sqrt{\lambda}}\cdot |\!|b|_{R^-\cup L^-}|\!|_{L^{\infty}} + |\!|b|_{R^-\cup L^-}|\!|^2_{L^{\infty}} \leq C_2\cdot K_N\cdot N^4\cdot \delta.\]

\vspace{1mm}

\noindent
(iii) Since $\nabla_A(\psi \cdot b) =  d\psi \otimes  b + \psi \cdot \nabla_Ab$ 
and $|d\psi| \leq 20N/\sqrt{\lambda}$,
\begin{eqnarray*}
|\!|\psi \cdot b|\!|_{L^p_1}&\leq& \frac{20N}{\sqrt{\lambda}}\cdot 
\mathrm{vol}(L^-\cup R^-)^{1/p} \cdot |\!|b|_{L^-\cup R^-}|\!|_{L^{\infty}} + |\!|b|\!|_{L^p_1} \\
&\leq& \mathrm{const}\cdot N^{1+4/p}\cdot \lambda^{-1/2 + 2/p}\cdot |\!|b|_{L^-\cup R^-}|\!|_{L^{\infty}} + |\!|b|\!|_{L^p_1}.
\end{eqnarray*}
From (\ref{implicit}) and (\ref{L^infty of b})
\[|\!|\psi \cdot b|\!|_{L^p_1} \leq C_3(N)\cdot \lambda^{2/p}\cdot \delta.\]

\noindent
(iv)\[F(A + a + \psi b) = F(A + a) + d_A(\psi b) + \psi[a \wedge b] + \psi^2 (b\wedge b).\]
Hence 
\[\norm{F(A + a + \psi b) - F(A + a)}_{L^2(U)} \leq 
\norm{\psi b}_{L^2_1} + \sqrt{2} \norm{a}_{L^4}\cdot \norm{b}_{L^4} + \sqrt{2} \norm{b}_{L^4}^2. \]
The arguments in (i) and (iii) show
\[ \norm{\psi b}_{L^2_1} \leq C(N) \cdot \lambda \cdot \delta, \quad \norm{b}_{L^4} \leq C'(N) \cdot \lambda \cdot \delta.\]
It follows that 
\[\norm{F(A + a + \psi b) - F(A + a)}_{L^2(U)} \leq C_4(N) \cdot \lambda \cdot \delta.\]
\end{proof}
Proposition \ref{prop: estimate for construction} is sufficient to prove Theorem \ref{thm: construction}, but Theorem \ref{thm: injectivity} needs more.
Next we will establish  a ``parametrized version" of Proposition \ref{prop: estimate for construction},
which is the key estimate for the proof of Theorem \ref{thm: injectivity}.
Here we follow the arguments in [D2, (4.24) Remarks].

Suppose that $a$ and $\sigma$ smoothly depend on a parameter $t \in [0,1]$: $a=a(t)$, $ \sigma = \sigma(t)$. (To be precise, it is enough to assume that $a : [0,1] \rightarrow L^{2p}$ and $ \sigma : [0,1] \rightarrow L^{\infty}$ are smooth maps.) 
Set $\delta = \delta(t) = |\!|\sigma(t) |\!|_{L^{\infty}(X)}$.
If $\delta \leq \kappa$ and $|\!|a|\!|_{L^{2p}(X)} \leq \eta$ for all $t$, 
($\kappa$ and $\eta$ are the constants in Proposition \ref{prop: estimate for construction}), 
then we have the following proposition.
\begin{Prop}
If we choose $\eta$ and $\lambda$ sufficiently small, 
then we have the following estimates:
\begin{eqnarray*}
&\mathrm{(i)}&
 \norm{\frac{\partial b}{\partial t}}_{L^{2p}(X)} \leq C_5(N)\cdot \lambda^{(p+2)/{2p}} \left( \norm{\frac{\partial \sigma}{\partial t}}_{L^{\infty}(X)} + \delta \cdot \norm{\frac{\partial a}{\partial t}}_{L^{2p}(X)} \right).\\
&\mathrm{(ii)}&
\norm{\frac{\partial \tau}{\partial t}}_{L^{\infty}(X)} \leq C_6 \cdot K_N \cdot N^4 \cdot \norm{\frac{\partial \sigma}{\partial t}}_{L^{\infty}(X)} + C_7(N)\cdot \lambda^{1/p} \cdot \delta \cdot \norm{\frac{\partial a}{\partial t}}_{L^{2p}(X)}.
\end{eqnarray*}
Here $C_5(N)$ and $C_7(N)$ are constants depending on $N$ and
$C_6$ is a constant independent of $N$.
\end{Prop}
\begin{proof}
(i) It is clear that $b : [0,1] \rightarrow L^p_1$ is smooth. 
Differentiating the equation $D_Ab + [a\wedge b]^+ +(b\wedge b)^+ = -\sigma$ with respect to $t$,
we get
\[D_A\left(\frac{\partial b}{\partial t}\right) + \left[\frac{\partial a}{\partial t} \wedge b\right]^+ +\left[a\wedge \frac{\partial b}{\partial t}\right]^+ +\left(\frac{\partial b}{\partial t}\wedge b\right)^+ +\left(b\wedge \frac{\partial b}{\partial t}\right)^+ = - \frac{\partial \sigma }{\partial t} .\]
Put $D'_A(\cdot) = D_A(\cdot) + [a\wedge \cdot]^+ +(\cdot \wedge b)^+ +(b\wedge \cdot)^+.$ 
 Then
\[D'_A\left(\frac{\partial b}{\partial t}\right) = -\left[\frac{\partial a}{\partial t} \wedge b\right]^+ -\frac{\partial \sigma}{\partial t}.\]
 Since $D_A$ is invertible and $\norm{a}_{L^{2p}}< \eta$ and $\norm{b}_{L^{2p}}< C_1(N)\cdot \lambda^{(2+p)/2p} \delta$ are sufficiently small, 
\begin{eqnarray*}
\norm{\frac{\partial b}{\partial t}}_{L^{2p}} &\leq& \mathrm{const} \cdot \norm{D'_A\left(\frac{\partial b}{\partial t}\right)}_{L^{4p/(p+2)}}\\
&\leq& \mathrm{const} \left( \norm{\frac{\partial \sigma }{\partial t}}_{L^{4p/(p+2)}} + \> \norm{\left[ \frac{\partial a}{\partial t} \wedge b\right]^+}_{L^{4p/(p+2)}}\right) \\
&\leq& C(N) \left(\lambda^{(2+p)/2p}\norm{\frac{\partial \sigma }{\partial t}}_{L^{\infty}} + \>
 \norm{\frac{\partial a}{\partial t}}_{L^{2p}} |\!|b|\!|_{L^{2p}}\right)\\
&\leq& C_5(N)\cdot \lambda^{(2+p)/2p}\left(\norm{\frac{\partial \sigma}{\partial t}}_{L^{\infty}} +  \>  \delta \cdot \norm{\frac{\partial a}{\partial t}}_{L^{2p}}\right).
\end{eqnarray*}
Here we have used the fact that $\partial \sigma /\partial t$ is supported in
$L^+\cup R^+$ whose volume is $O(\lambda^2)$.

\vspace{1mm}
\noindent
(ii) Differentiating $\tau = (d\psi \wedge b)^+ + \psi(\psi -1)(b\wedge b)^+$ with respect to $t$,
we have 
\[ \frac{\partial \tau}{\partial t} = \left(d\psi \wedge \frac{\partial b}{\partial t}\right)^+ +\psi(\psi -1)\{\left(\frac{\partial b}{\partial t} \wedge b\right)^+ + \left(b \wedge \frac{\partial b}{\partial t}\right)^+\}.\]
From (\ref{cutoff}) and (\ref{L^infty of b})
\begin{equation}\label{estimate of tau/t}
 \begin{split}
 \norm{\frac{\partial \tau}{\partial t}}_{L^{\infty}} &\leq \> \frac{20N}{\sqrt{\lambda}}
 \norm{\left.\frac{\partial b}{\partial t}\right|_{L^- \cup R^-}}_{L^{\infty}}
 + 2\sqrt{2}\cdot \norm{b|_{L^- \cup R^-}}_{L^{\infty}}
 \norm{\left.\frac{\partial b}{\partial t}\right|_{L^- \cup R^-}}_{L^{\infty}} \\
 & \leq \> \frac{20N}{\sqrt{\lambda}}
 \norm{\left.\frac{\partial b}{\partial t}\right|_{L^- \cup R^-}}_{L^{  \infty}}
 + C(N)\cdot \sqrt{\lambda}\cdot \delta \cdot 
 \norm{\left.\frac{\partial b}{\partial t}\right|_{L^- \cup R^-}}_{L^{\infty}}.
 \end{split}
\end{equation}
Set $b = \tilde{b} + \beta$ with $D_A \tilde{b} = - \sigma$ 
as in the proof of Proposition 3.1 (ii).
Then
\begin{eqnarray*}
\dfrac{\partial \tilde{b}}{\partial t} (x) &=& -\int_{X} G(x,y)\cdot \dfrac{\partial \sigma}{\partial t} (y)\, dy, \quad (G(x,y) \text{ is the Green kernel}),\\
&=& - \int_{L^+} G(x,y)\cdot \dfrac{\partial \sigma}{\partial t} (y)\, dy 
-\int_{R^+}  G(x,y)\cdot \dfrac{\partial \sigma}{\partial t} (y) \, dy.
\end{eqnarray*}
This can be estimated in the same way as before and we get
\[\norm{\left.\frac{\partial \tilde{b}}{\partial t}\right|_{L^- \cup R^-}}_{L^{\infty}}\leq \mathrm{const}\cdot K_N \cdot N^3 \cdot \sqrt{\lambda}\cdot \norm{\frac{\partial \sigma}{ \partial t}}_{L^{\infty}}.\]
We have $\partial b/\partial t = \partial \tilde{b}/\partial t + \partial \beta/\partial t$ 
and $\partial \beta/\partial t$ satisfies 
\[D_A\left(\frac{\partial \beta}{\partial t}\right) + \left[\frac{\partial a}{\partial t} \wedge b\right]^+ + \left[a\wedge \frac{\partial b}{\partial t}\right]^+ +\left(\frac{\partial b}{\partial t} \wedge b\right)^+ +\left(b\wedge \frac{\partial b}{\partial t}\right)^+ = 0.\]
Using the above (i) and the Sobolev embedding, we get 
\begin{eqnarray*}
\norm{\frac{\partial \beta}{\partial t}}_{L^{\infty}} &\leq& \mathrm{const}\cdot \norm{D_A\left(\frac{\partial \beta}{\partial t}\right)}_{L^6}, \quad (L^6_1\hookrightarrow L^{\infty}),\\
&\leq& C(N)\cdot \lambda^{(2+p)/2p}\left(\norm{\frac{\partial \sigma}{\partial t}}_{L^{\infty}}+ \delta \cdot \norm{\frac{\partial a}{\partial t}}_{L^{2p}}\right). 
\end{eqnarray*}
Hence 
\[\norm{\left.\frac{\partial b}{\partial t}\right|_{L^- \cup R^-}}_{L^{\infty}} \leq C\cdot K_N  N^3 \sqrt{\lambda}\norm{\frac{\partial \sigma }{\partial t}}_{L^{\infty}} + C(N)\cdot \lambda^{(2+p)/2p} \cdot \delta \norm{\frac{\partial a}{\partial t}}_{L^{2p}}.\]
Substituting this into (\ref{estimate of tau/t}), we conclude that 
\[\norm{\frac{\partial \tau}{\partial t}}_{L^{\infty}} \leq C_6 \cdot K_N N^4 \norm{\frac{\partial \sigma}{\partial t}}_{L^{\infty}} + C_7(N)\cdot \lambda^{1/p}\cdot \delta  \norm{\frac{\partial a}{\partial t}}_{L^{2p}}.\]
\end{proof}
Here we fix the value of N. Any value of N with
\[C_2\cdot K_N \leq \frac{1}{2}, \quad C_6\cdot K_N \leq \frac{1}{2},\]
will do. (Note that $K_N=N/(N-N^{-1})^3$ converges to zero as $N \rightarrow \infty$.)
So let N be a large positive number which satisfies the above.
(This choice depends on the gluing data $(X, x^L, x^R, E, A)$.) 
Then Proposition 3.1 (ii) and Proposition 3.2 (ii) become:
\begin{eqnarray*}
&\text{(ii)$'$}&
 \norm{\tau}_{L^{\infty}} \leq \frac{1}{2}\cdot N^4 \cdot \delta ,\\
&\text{(ii)$'$}&
 \norm{\frac{\partial \tau}{\partial t}}_{L^{\infty}} \leq \frac{1}{2}\cdot N^4 \cdot \norm{\frac{\partial \sigma}{\partial t}}_{L^{\infty}} + C_7(N)\cdot \lambda^{1/p} \cdot \delta \cdot
 \norm{\frac{\partial a}{\partial t}}_{L^{2p}}.
\end{eqnarray*}
\begin{Rem}
In fact, the parameter $N$ can be chosen independent of the gluing data $(X, x^L, x^R, E, A)$. 
See [D2, p. 308, (4.14)] and the proof of [D2, (4.22) Lemma].
But we don't need this fact in this paper.
\end{Rem}


\section{Donaldson's alternating method}
In this section we prove Theorem \ref{thm: construction}. Let $(X_i, x^L_i, x^R_i, E_i, A_i)_{i \in \mathbb{Z}}$ be a sequence of gluing data of finite type as in Section 2.
The finite type condition means that 
we can apply the estimates in Section 3 to all $(X_i, x^L_i, x^R_i, E_i, A_i)$ 
with the \textit{same} constants. 
(This is the reason for introducing this notion.)

To begin with, 
we make one remark on the conformal transformation (\ref{conformal map}) 
which identifies $\Omega^R_i$ and $\Omega^L_{i+1}$:
\begin{Rem}\label{rem: conformal map}
The coordinate transformation (\ref{conformal map}) is only conformal and not an isometry. 
Hence the metric of $X_i$ is different from the metric of $X_{i+1}$ 
over the identified region $\Omega^R_i \cong \Omega^L_{i+1}$. 
The derivative of (\ref{conformal map}) over 
$\Omega^R_i = \{kN^{-1}\sqrt{\lambda} < |\xi| < k^{-1}N\sqrt{\lambda}\}$ 
has the following bound independent of $\lambda$:
\[\norm{\frac{\partial\eta}{\partial \xi}}= \frac{\lambda}{\abs{\xi}^2}\leq \frac{N^2}{k^2},\quad 
 \norm{\frac{\partial \xi}{\partial \eta}} = \frac{\lambda}{\abs{\eta}^2} \leq
\frac{N^2}{k^2}.
\]
Here $\norm{\cdot}$ is the operator norm.
In particular, if $\xi \in R^-_i = L^+_{i+1}$, then
\[ \norm{\frac{\partial \xi}{\partial \eta}} \leq \frac{1}{N^2}.\]
It follows that
\[ \norm{\omega}_{L^{\infty}(L^+_{i+1})}\leq \frac{1}{N^4}\norm{\omega}_{L^{\infty}(R^-_i)}\]
 for a 2-form $\omega$ under the identification $R^-_{i} = L^+_{i+1}$.
\end{Rem}
Let $\rho=(\rho_i)_{i\in \mathbb{Z}} \in \mathrm{GlP}$ be a gluing parameter, 
($\rho_i \in \mathrm{Hom}_{G}((E_i)_{x^R_i},(E_{i+1})_{x^L_{i+1}})$),
and $\bm{E}(\rho)$ be the principal $G$ bundle over $\XX$ defined in 
Section \ref{subsection: gluing construction}. 
As we explained in Section \ref{subsection: outline},
we will inductively construct a sequence of connections 
$\bm{A}^{(0)}(\rho)$, $\bm{A}^{(1)}(\rho)$, $\bm{A}^{(2)}(\rho)$, $ \cdots$, on $\bm{E}(\rho)$ which converges to an ASD connection $\bm{A}(\rho)$. 
This is a generalization of the alternating method in [D2, IV (iv)].
At the $n$-th stage of the iteration we will have a connection $\bm{A}^{(n)}(\rho)$ which is expressed by 
\[\bm{A}^{(n)}(\rho) = (A_i + a^{(n)}_i(\rho))_{i \in \mathbb{Z}}.\]
Here $a^{(n)}_i(\rho)$ is an element of $\Omega^1(\mathrm{ad}E_i)$ over $U_i$ and we have
$A_i + a^{(n)}_i(\rho) = A_{i+1}+a^{(n)}_{i+1}(\rho)$ over $U_i\cap U_{i+1}$ 
under the identification $\rho_i$.
Then they compatibly define the connection $\bm{A}^{(n)}(\rho)$ on $\bm{E}(\rho)$.


Let us consider the following ``inductive hypotheses":
\begin{Hyp}\label{inductive hypotheses}
(i) Let $\sigma^{(n)} := F^{+}(\bm{A}^{(n)}(\rho))$ be the self-dual curvature of $\bm{A}^{n}(\rho)$, 
and set 
\[\delta_n := \begin{cases}
\mathrm{sup}_{i \in \mathbb{Z}} |\!|\sigma^{(n)}_{2i}|\!|_{L^{\infty}(X_{2i})} & \text{if }n \text{ is even.}\\
\mathrm{sup}_{i \in \mathbb{Z}}|\!|\sigma^{(n)}_{2i+1}|\!|_{L^{\infty}(X_{2i+1})} & \text{if }n \text{ is odd.}
\end{cases}
\]
Here $\sigma^{(n)}_{i}$ is the restriction of $\sigma^{(n)}$ on $U_i$. 
(We extend $\sigma^{(n)}_i$ to $X_i$ by zero.) 
We suppose that $\sigma^{(n)}$ is supported as follows:  
\[\mathrm{supp}(\sigma^{(n)}) \subset \begin{cases}
\coprod_{i\in \mathbb{Z}}(L^+_{2i}\cup R^{+}_{2i}) & \text{if }n\text{ is even.}\\
\coprod_{i\in \mathbb{Z}}(L^+_{2i+1}\cup R^{+}_{2i+1}) & \text{if }n\text{ is odd.}
\end{cases}
\]
And we suppose that $\delta_n$ has the following bound:
\[ \delta_n \leq \kappa , \quad (\kappa \text{ is the constant in Proposition 3.1 and will be given later}).\]

\noindent
(ii) We suppose that $a^{(n)}_i$ has the following bound:
\[ |\!|a^{(n)}_i|\!|_{L^{2p}(U_i)} \leq \eta , \quad 
(\eta \text{ is also the constant in Proposition 3.1}) .\]
\end{Hyp}
If these conditions hold at the $n$-th stage of the iteration and $\lambda \ll 1$, 
then we can pass to the $(n+1)$-th stage as follows.
If $n$ is even, then we use the construction in Section 3 on each $X_{2i}$ and get $b^{(n)}_{2i}$. 
We define $a^{(n)}_{i}(\rho)$ over $U_i$ by
\begin{equation}\label{iteration1}
 \begin{split}
   &a^{(n+1)}_{2i}(\rho) :=  a^{(n)}_{2i}(\rho) + \psi_{2i}\cdot b^{(n)}_{2i}, \\  
   &a^{(n+1)}_{2i+1}(\rho) :=  a^{(n)}_{2i+1}(\rho) + \psi_{2i}\cdot \rho_{2i}b^{(n)}_{2i}
\rho_{2i}^{-1} + \psi_{2i+2}\cdot \rho_{2i+1}^{-1}b^{(n)}_{2i+2}\rho_{2i+1}. 
 \end{split}
\end{equation}
Here $\psi_i$ is a cut-off over $X_i$ as in Section 3. Note that $\psi_{2i}\cdot \rho_{2i}b^{(n)}_{2i}\rho_{2i}^{-1}$ and $\psi_{2i+2}\cdot \rho_{2i+1}^{-1}b^{(n)}_{2i+2}\rho_{2i+1}$ are well-defined over $U_{2i+1}$.
From the arguments in Section 3 and Proposition \ref{prop: estimate for construction} (ii)$'$ at the end of Section 3, we have
\[
\mathrm{supp}(\sigma^{(n+1)}) \subset \coprod_{i \in \mathbb{Z}}(R^-_{2i}\cup L^-_{2i}) 
= \coprod_{i \in \mathbb{Z}}(L^+_{2i+1}\cup R^+_{2i+1}), \]
\[ |\!|\sigma^{(n+1)}_{2i}|\!|_{L^{\infty}(X_{2i})} \leq \frac{1}{2}\cdot N^4 \delta_n. \]
Taking account of Remark \ref{rem: conformal map}, we have 
\[|\!|\sigma^{(n+1)}_{2i+1}|\!|_{L^{\infty}(X_{2i+1})} \leq \frac{1}{2}\delta_n.\]
Thus
\begin{equation}
\delta_{n+1} \leq \frac{1}{2}\delta_n. \label{er.term}
\end{equation}
If $n$ is odd, then we use Proposition \ref{prop: estimate for construction} over each $X_{2i+1}$
and construct a similar procedure to pass to the $(n+1)$-th stage.
The estimate (\ref{er.term}) holds as long as we can continue this iteration. 
Hence
\begin{equation}
\delta_n \leq \frac{1}{2^n}\delta_0. \label{geo.decay}
\end{equation} 

We start the iteration by defining $\bm{A}^{(0)}(\rho) = (A_i + a^{(0)}_i(\rho))_{i \in \mathbb{Z}}$ by
\begin{equation}\label{ini.value1}
 \begin{split}
 &a^{(0)}_{2i+1}(\rho) := 
 (\psi_{2i+1}-1)(A_{2i+1} - \rho_{2i}A_{2i}\rho_{2i}^{-1} -\rho_{2i+1}^{-1}A_{2i+2}\rho_{2i+1}), \\
 &a^{(0)}_{2i}(\rho) := 
 \psi_{2i+1}(-A_{2i} + \rho_{2i}^{-1}A_{2i+1}\rho_{2i}) + 
 \psi_{2i-1}(-A_{2i} + \rho_{2i-1}A_{2i-1}\rho_{2i-1}^{-1}). 
 \end{split}
\end{equation}
Here each $A_i$ is the connection matrix in the exponential gauge centered at $x^L_i$ and $x^R_i$.
These expressions are well-defined on $U_{2i+1}$ or $U_{2i}$ and compatible over the overlap regions.
The self-dual curvature of $\bm{A}^{(0)}(\rho)$ is supported in $\coprod (L^-_{2i+1}\cup R^-_{2i+1}) = \coprod (L^+_{2i}\cup R^+_{2i})$. 
Since each $\abs{A_i}$ is $O(\sqrt{\lambda})$ and 
$\abs{d\psi_i}$ is $O(1/\sqrt{\lambda})$ (cf. (\ref{cutoff})), 
we have a constant $\kappa$ independent of $\lambda$ such that 
\begin{equation}
 |\!|\sigma^{(0)}_{i}|\!|_{L^{\infty}(X_{i})} \leq \kappa. \label{ini.error}
\end{equation}
We take this $\kappa$ as the constant in Proposition \ref{prop: estimate for construction}.
The condition (i) in Hypotheses \ref{inductive hypotheses} is satisfied at $n=0$.
Since $a^{(0)}_i$ is supported in the overlap whose volume is $O(\lambda^2)$, we have
\begin{equation}
|\!|a^{(0)}_i|\!|_{L^{2p}(U_i)} \leq \mathrm{const}\cdot \lambda^{(2+p)/2p}. \label{ini.estimate1}
\end{equation}
Therefore the condition (ii) is also satisfied at $n=0$ if we choose $\lambda \ll 1$. 
Then we can start the iteration.

The geometric decay (\ref{geo.decay}) shows that the condition (i) in Hypotheses 4.2 is always satisfied as long as the condition (ii) is valid.
At the $n$-th stage of the iteration, if $n$ is even, Proposition 3.1 (i) shows \[|\!|a^{(n+1)}_{2i} - a^{(n)}_{2i}|\!|_{L^{2p}(U_{2i})} \leq \mathrm{const}\cdot \lambda^{(2+p)/2p}\delta_n \leq \mathrm{const}\cdot \lambda^{(2+p)/2p}\frac{\delta_0}{2^n}.\]
From Remark \ref{rem: conformal map},
\[ |\!|a^{(n+1)}_{2i+1} - a^{(n)}_{2i+1}|\!|_{L^{2p}(U_{2i+1})} \leq \mathrm{const} \cdot \lambda^{(2+p)/2p}\frac{\delta_0}{2^n}.\]
We have similar estimates when $n$ is odd, and
\begin{equation}\label{total error}
 \mathrm{const}\cdot \lambda^{(2+p)/2p} + \mathrm{const}\cdot \lambda^{(2+p)/2p} \sum_{n\geq 0}
 \frac{\delta_0}{2^n} \leq \mathrm{const} \cdot \lambda^{(2+p)/2p}\ll 1. 
\end{equation}
Hence the condition (ii) in Hypotheses 4.2 is satisfied for any $n$ if we choose $\lambda$ small enough. Thus we can continue the iteration indefinitely and 
the sequence $\{ a^{(n)}_i(\rho) \}_{n=0}^{\infty}$ has a limit $ a_i(\rho)$ in $L^{2p}(U_i)$.
 
Proposition 3.1 (iii) and an argument similar to the above show that $\{ a^{(n)}_i(\rho) \}_{n=0}^{\infty}$ converges to $a_i(\rho)$ also in $L^p_1(U_i)$. 
From the geometric decay (\ref{geo.decay}), the self-dual curvature of $A_i + a_i(\rho)$ becomes $0$.
Thus we get an ASD connection
\[\bm{A}(\rho) = (A_i + a_i(\rho))_{i\in \mathbb{Z}}.\]
(\ref{total error}) shows that this is close to $A_i$ over each $U_i$:
\begin{equation}
|\!|a_i|\!|_{L^{2p}(U_i)} \leq \mathrm{const}\cdot \lambda^{(2+p)/2p}. \label{finalterm}
\end{equation}
 So we get a conclusion: 
\begin{Thm}
For all sufficiently small $\lambda$, $\bm{A}(\rho)$ is an ASD connection on $\bm{E}(\rho)$ such that 
\[|\!|\bm{A}(\rho) - A_i|\!|_{L^{2p}(U_i)} \leq C_8 \cdot \lambda^{(2+p)/2p}
\quad \text{for all } i \in \mathbb{Z}.\] 
Here $C_8$ is a constant independent of $i$, $\rho$, $\lambda$.
\end{Thm}
In general this ASD connection $\bm{A}(\rho)$ has infinite energy:
\begin{Thm}
If there are infinitely many non-flat ASD connections $A_i$ in the given sequence of gluing data, 
then the $L^2$-energy of $\bm{A}(\rho)$ is infinite: $\norm{F(\bm{A}(\rho))}_{L^2}= + \infty$.
\end{Thm}
\begin{proof}
We will show that the loss of energy over each $U_i$ during the iteration is sufficiently small. 
At the $n$-th stage of the iteration, if $n$ is even, we can apply Proposition 3.1 (iv) to (\ref{iteration1}) and get 
\[ |\!|F(A_{2i} + a^{(n+1)}_{2i}) - F(A_{2i} + a^{(n)}_{2i})|\!|_{L^2(U_{2i})} \leq \mathrm{const}\cdot \lambda \cdot \delta_n \leq \mathrm{const} \cdot \frac{\lambda}{2^n}.\]
Since $F(A_{2i+1} + a^{(n+1)}_{2i+1}) - F(A_{2i+1} + a^{(n)}_{2i+1})$ is supported in the overlaps 
and the $L^2$-norm of a 2-form is conformally invariant,
\[|\!|F(A_{2i+1} + a^{(n+1)}_{2i+1}) - F(A_{2i+1} + a^{(n)}_{2i+1})|\!|_{L^2(U_{2i+1})} \leq \mathrm{const}\cdot \frac{\lambda}{2^n}.\]
We have the same estimates in the case that $n$ is odd. Hence
\[|\!| F(A_{i} + a_{i}) - F(A_{i} + a^{(0)}_{i})|\!|_{L^2(U_{i})} \leq \mathrm{const} \cdot \lambda \sum_{n\geq 0}\frac{1}{2^n} \leq \mathrm{const} \cdot \lambda. \]
Next we will show that $|\!| F(A_{i} + a^{(0)}_{i}) - F(A_{i})|\!|_{L^2(U_{i})}$ is also small. 
An estimate similar to (\ref{ini.error}) shows 
\[|F(A_i + a_i^{(0)}) - F(A_i)| \leq \mathrm{const}, \quad (\text{independent of }\lambda).\] 
Since this is supported in the overlaps, we get
\[|\!|F(A_i + a_i^{(0)}) - F(A_i)|\!|_{L^2(U_i)} \leq \mathrm{const} \cdot \lambda.\]
Hence
\[|\!| F(A_{i} + a_{i}) - F(A_{i})|\!|_{L^2(U_{i})} \leq \mathrm{const} \cdot \lambda.\]

From the finite type condition, we have a positive constant $M$ independent of $\lambda$ such that $|\!|F(A_i)|\!|_{L^2(U_i)} \geq M$ if $A_i$ is not flat and $\lambda$ is sufficiently small.
Therefore if we take $\lambda \ll 1$, then
\[|\!|F(A_i + a_i)|\!|_{L^2(U_i)} \geq \frac{1}{2}M \]
for non-flat $A_i$.
 Since there are infinitely many non-flat $A_i$, we conclude that
\[|\!|F(\bm{A}(\rho))|\!|_{L^2} = + \infty.\]  
\end{proof}


\section{Moduli Problem}
We will prove Theorem \ref{thm: injectivity} in this section. The proof is a generalization of [D2, (4.31) Lemma]. We always suppose that $\lambda$ is sufficiently small.
 
We need the following elementary result.
\begin{Lem}\label{lem:recurrence}
Let $K$ be a positive number and $\varepsilon$ be a positive number with $\varepsilon \leq 1/100$.
Let $\{ \alpha_n\}, \{ s_n\}$, $(n=0,1,2,\cdots)$, be sequences of non-negative numbers which satisfy 
the following recurrence inequalities:
\begin{equation*}
 \begin{split}
 \alpha_{n+1} - \alpha_n &\leq \varepsilon(s_n + \frac{1}{2^n}\alpha_n +\frac{1}{2^n} K),  \\
 s_{n+1} &\leq \frac{1}{2}s_n + \frac{1}{2^n} \alpha_n + \frac{1}{2^{n+1}} K. 
 \end{split}
\end{equation*}
Suppose that 
\[\alpha_0 \leq \varepsilon \cdot K , \quad s_0 \leq K.\]
Then we get 
\[ \alpha_n \leq 10 \varepsilon \cdot K \quad \text{for all } n.\]
\end{Lem}
\begin{proof}
Set $t_n := 2^n s_n$. Then we have $t_0 \leq K$, and the above recurrence inequalities become 
\begin{equation*}
 \begin{split}
 \alpha_{n+1} - \alpha_n &\leq \frac{\varepsilon}{2^n}(t_n + \alpha_n + K), \\
 t_{n+1} &\leq t_n + 2 \alpha_n + K.
 \end{split}
\end{equation*}
We will inductively prove $\alpha_n \leq 10\varepsilon \cdot K$ for all $n$. 
Suppose that this is true if $n \leq n_0$ for some $n_0$.
Then 
\[t_n \leq K + 2 (10 \varepsilon K)n + K n  =
  K \{ 1 + (1 + 20 \varepsilon)n \},  \quad(n \leq n_0 + 1).  \]
Hence 
\begin{eqnarray*}
\alpha_{n_0 +1} - \alpha_0 &=& \sum_{n=0}^{n_0} (\alpha_{n+1} - \alpha_{n})\\
     &\leq& \sum_{n\geq 0} \frac{\varepsilon}{2^n}\left[K \{ 1 + (1 + 20 \varepsilon)n \} + 10 \varepsilon K + K \right] \\
     &=& \varepsilon K (6 + 60\varepsilon) \leq 7 \varepsilon K.  
\end{eqnarray*}
Thus 
\[ \alpha_{n_0 + 1} \leq \varepsilon K + 7 \varepsilon K \leq 10 \varepsilon K. \]
\end{proof}


\begin{Prop}
For two \textit{gluing parameters} $\rho = (\rho_i)_{i\in\mathbb{Z}}$ and $\rho' = (\rho'_i)_{i\in \mathbb{Z}}$, let $K:= \mathrm{sup}_{i\in \mathbb{Z}}|\rho_i - \rho'_i|$. Then we have
\[ |\!|a_i(\rho) - a_i(\rho')|\!|_{L^{2p}(U_i)} \leq C_9 \cdot \lambda^{(2+p)/2p} K
\quad \text{for all } i \in \mathbb{Z}. \]
Here $C_9$ is a constant independent of $i$, $\lambda$, $\rho$, $\rho'$.
\end{Prop}
\begin{proof}
Let $\bar{\rho}(t) = (\bar{\rho_i}(t))_{i\in\mathbb{Z}}$, $(0\leq t \leq 1)$, be a path in $\mathrm{GlP}$ such that each $\bar{\rho_i}(t)$ is the geodesic from $\rho_i$ to $\rho'_i$ in $G$.
Since the distance $d(\rho_i, \rho'_i)$ is equivalent to $|\rho_i - \rho'_i|$:
\[ C^{-1} |\rho_i - \rho'_i| \leq d(\rho_i, \rho'_i) \leq C |\rho_i - \rho'_i| \quad \text{for some constant } C, \]
we have $|d \bar{\rho_i}/dt| = d(\rho_i, \rho'_i) \leq C\cdot K$.

Let $\bm{A}^{(n)}(\bar{\rho}(t)) = ( A_i + a_i^{(n)}(\bar{\rho}(t)))_{i \in \mathbb{Z}}$ 
be the connection on $\bm{E}(\bar{\rho}(t))$ constructed by the iterative process in Section 4; this is parametrized by $t$ and we can apply Proposition 3.2. 
We define $\alpha_n = \alpha_n(t)$ and $s_n = s_n(t)$ by
\[ \alpha_n := \mathrm{sup}_{i \in \mathbb{Z}} \norm{\deldel a^{(n)}_i}_{L^{2p}(X_i)},\quad
   s_n := \begin{cases}
          \mathrm{sup}_{i \in \mathbb{Z}} \norm{\deldel \sigma_{2i}^{(n)}}_{L^{\infty}(X_{2i})}          & \text{ if } n \text{ is even.} \\
          \mathrm{sup}_{i \in \mathbb{Z}} \norm{\deldel \sigma_{2i + 1}^{(n)}}_{L^{\infty}(X_{2i+1})}  & \text{ if } n \text{ is odd.}
          \end{cases}
\]
From the definition of $a_i^{(0)}$ in (\ref{ini.value1}), 
\[ \abs{\deldel a^{(0)}_i} \leq \mathrm{const}\cdot \sqrt{\lambda} \cdot K. \]
Here we used the fact that the connection matrices in the exponential gauges are $O(\sqrt{\lambda})$.
Since $\partial a^{(0)}_i / \partial t$ is supported in the overlaps, we get
\[ \norm{\deldel a^{(0)}_i}_{L^{2p}(X_i)}  \leq \mathrm{const} \cdot \lambda^{(2+p)/2p} K.\]
In a similar way, we have 
\[ \norm{\deldel \sigma^{(0)}_i}_{L^{\infty}(X_i)} \leq \mathrm{const} \cdot K.\]
Hence
\begin{equation}
 \alpha_0 \leq \mathrm{const} \cdot \lambda^{(2+p)/2p} K, \quad s_0 \leq \mathrm{const} \cdot K. 
\label{ini.value3}
\end{equation}

Next we will establish recurrence inequalities on $\alpha_n$ and $s_n$.
If $n$ is even, then Proposition 3.1 and 3.2 give bounds for $|\!| b^{(n)}_{2i} |\!|_{L^{2p}(X_{2i})}$ and $|\!| \partial b^{(n)}_{2i}/ \partial t |\!|_{L^{2p}(X_{2i})}$. 
Then, from the inductive definitions (\ref{iteration1}) and the geometric decay (\ref{geo.decay}),
\begin{equation*}
 \begin{split}
 &\norm{\deldel a^{(n+1)}_{2i}}_{L^{2p}(X_{2i})} \leq \alpha_n + \mathrm{const} \cdot 
 \lambda^  {(2+p)/2p}(s_n + \frac{1}{2^n} \alpha_n), \\
 &\norm{\deldel a^{(n+1)}_{2i+1}}_{L^{2p}(X_{2i+1})} \leq \alpha_n + \mathrm{const} \cdot 
 \frac{\lambda^{(2+p)/2p}}{2^n} K + \mathrm{const} \cdot \lambda^{(2+p)/2p} 
 (s_n + \frac{1}{2^n} \alpha_n).
 \end{split}
\end{equation*}
It follows that 
\begin{equation}
\alpha_{n+1} - \alpha_{n} \leq \mathrm{const} \cdot \lambda^{(2+p)/2p}(s_n + \frac{1}{2^n} \alpha_n + \frac{K}{2^n}).   \label{iteration3}
\end{equation}
Similarly, from Proposition 3.2 (ii)$'$ at the end of Section 3,
\[\norm{\deldel \sigma^{(n+1)}_{2i}}_{L^{\infty}(X_{2i})} \leq \frac{N^4}{2} s_n + \mathrm{const} \cdot \frac{\lambda^{1/p}}{2^n} \alpha_n. \]
Since $\sigma^{(n+1)}_{2i+1} = \bar{\rho}_{2i} \sigma^{(n+1)}_{2i} \bar{\rho}_{2i}^{-1}$ on 
$U_{2i} \cap U_{2i+1}$ and 
$\sigma^{(n+1)}_{2i+1} = \bar{\rho}_{2i+1}^{-1} \sigma^{(n+1)}_{2i+2} \bar{\rho}_{2i+1}$
on $U_{2i+1}\cap U_{2i+2}$, we get
\begin{eqnarray*} 
\norm{\deldel \sigma^{(n+1)}_{2i+1}}_{L^{\infty}(X_{2i+1})} &\leq& \mathrm{const} \cdot \delta_{n+1} \cdot K + \frac{1}{2} s_n + \mathrm{const} \cdot \frac{\lambda^{1/p}}{2^n} \alpha_n \\
&\leq& \mathrm{const} \cdot \frac{K}{2^{n+1}}  + \frac{1}{2} s_n + \mathrm{const} \cdot \frac{\lambda^{1/p}}{2^n} \alpha_n.
\end{eqnarray*}
Here we have used Remark \ref{rem: conformal map}. 
This shows
\begin{equation}
s_{n+1} \leq \frac{1}{2} s_n + \mathrm{const} \cdot \frac{\lambda^{1/p}}{2^n} \alpha_n + \mathrm{const} \cdot \frac{K}{2^{n+1}}.  \label{iteration4}
\end{equation}
The inequalities (\ref{iteration3}), (\ref{iteration4}) hold also in the case that $n$ is odd.
So we can apply Lemma \ref{lem:recurrence}
with initial values (\ref{ini.value3}) and $\varepsilon = \mathrm{const} \cdot \lambda^{(2+p)/2p}$. 
Then we get
\[ \alpha_n \leq \mathrm{const} \cdot \lambda^{(2+p)/2p} K. \]
Therefore
\[ |\!| a^{(n)}_i(\rho) - a^{(n)}_i(\rho')|\!|_{L^{2p}(U_i)} = \norm{\int_{0}^{1} \deldel a^{(n)}_i dt}_{L^{2p}(U_i)} \leq \mathrm{const} \cdot \lambda^{(2+p)/2p} K,\]
for all $i$ and $n$.
So we get a conclusion:
\[ |\!|a_i(\rho) - a_i(\rho')|\!|_{L^{2p}(U_i)} \leq \mathrm{const} \cdot \lambda^{(2+p)/2p} K. \]
\end{proof}

The following is Theorem \ref{thm: injectivity}.
\begin{Thm}
For two \textit{gluing parameters} $\rho$ and $\rho'$, $\bm{A}(\rho)$ is gauge equivalent to $\bm{A}(\rho')$ if and only if $[\rho] = [\rho']$ in $\mathrm{EGlP}$.
\end{Thm}
\begin{proof}
If $[\rho] = [\rho']$ in $\mathrm{EGlP}$, we have $\varepsilon_i \in C(G)$ 
with $\rho'_i = \varepsilon_i \rho_i$ for all $i\in \mathbb{Z}$.
Since the center $C(G)$ acts trivially on connections, 
it is clear from the construction in Section 4 that 
we have $a_i(\rho) = a_i(\rho')$.
And we can choose $\gamma_i \in C(G)$ such that $\gamma_{i+1} = \varepsilon_i \gamma_i$.
Then these $\gamma_i$ give a gauge transformation 
$h = (\gamma_i)_{i\in \mathbb{Z}} :\bm{E}(\rho) \rightarrow \bm{E}(\rho')$ with 
$h. \bm{A}(\rho) = \bm{A}(\rho')$.

Conversely, suppose that there is a gauge transformation $ g: \bm{E}(\rho) \rightarrow \bm{E}(\rho')$ 
with $ g. \bm{A}(\rho) = \bm{A}(\rho')$.
This means that there is a gauge transformation $g_i: E_i \rightarrow E_i$ over each $U_i$ such that
\begin{equation*}
 \begin{split}
 &d_{A_i}{g_i} = g_i a_i(\rho) - a_i(\rho' )g_i, \\
 &g_{i+1}\rho_i = \rho'_i g_i.
 \end{split}
\end{equation*} 
Because $A_i$ is irreducible, we have a constant C independent of $\lambda$ such that 
\[\min_{\gamma \in C(G)}|\!|\gamma g_i - 1|\!|_{C^0(U_i)}
\leq C \cdot |\!|d_{A_i} g_i|\!|_{L^{2p}(U_i)}. \] 
The proof of this inequality is standard. For example, see  [D2, (4.31) Lemma].
We can suppose that $|\!|g_i - 1|\!|_{C^0(U_i)} \leq C \cdot |\!|d_{A_i} g_i|\!|_{L^{2p}(U_i)}$ 
without loss of generality. 
From (\ref{finalterm}) and Proposition 5.2,
\begin{eqnarray*}
|\!|g_i - 1|\!|_{C^0} &\leq& C\cdot |\!|g_i a_i(\rho) - a_i(\rho')g_i|\!|_{L^{2p}} \\
&\leq& C\cdot |\!|g_i - 1|\!|_{C^0}(|\!|a_i(\rho)|\!|_{L^{2p}} + |\!|a_i(\rho')|\!|_{L^{2p}}) + C \cdot |\!|a_i(\rho) - a_i(\rho')|\!|_{L^{2p}} \\
&\leq& \mathrm{const} \cdot\lambda^{(2+p)/2p} |\!|g_i - 1|\!|_{C^0} + \mathrm{const}\cdot
\lambda^{(2+p)/2p} K.
\end{eqnarray*}
Since $\lambda$ is sufficiently small, we get
\[|\!|g_i - 1|\!|_{C^0} \leq \mathrm{const} \cdot \lambda^{(2+p)/2p} K.\]
On the other hand, the ``compatibility condition" $g_{i+1}\rho_i = \rho'_i g_i$ means 
\[\rho'_i - \rho_i = (g_{i+1}-1)\rho_i g^{-1}_i + \rho_i g^{-1}_i (1-g_i). \]
Hence 
\[ |\rho'_i - \rho_i| \leq \mathrm{const} \cdot (|\!|g_{i+1} -1|\!|_{C^0(U_{i+1})} + |\!|g_i - 1 |\!|_{C^0(U_{i})})  \leq \mathrm{const} \cdot \lambda^{(2+p)/2p} K. \]
Since this holds for all $i$, it follows 
\[K \leq \mathrm{const} \cdot \lambda^{(2+p)/2p} K.\]
If we choose $\lambda$ small enough, then this shows $K=0$ and $\rho = \rho'$.
\end{proof}

\begin{Rem}
So far we have studied only acyclic connections. 
But this is just for simplicity. 
The condition (iv) in Definition \ref{def: gluing data} can be replaced with
\[ \text{(iv)' } A \text{ is a regular ASD connection on } E. \]
If we glue an infinite number of regular ASD connections, 
then the resulting ASD connections are parametrized by
\[ \mathrm{EGlP} \times \prod_{i\in \mathbb{Z}}B_i
\subset \mathrm{EGlP} \times \prod_{i\in \mathbb{Z}}H^1_{A_i}.\]
Here $B_i$ is a small ball centered at the origin in $H^1_{A_i} := \mathrm{ker} D_{A_i}$. This generalization is a simple application of the methods in this paper and [D2, pp. 324-325]. 
The proof is routine, so I omit the details.
\end{Rem}


\vspace{10mm}

\address{ Masaki Tsukamoto \endgraf
Department of Mathematics, Faculty of Science \endgraf
Kyoto University \endgraf
Kyoto 606-8502 \endgraf
Japan
}

\textit{E-mail address}: \texttt{tukamoto@math.kyoto-u.ac.jp}

\end{document}